# CRITERIA OF PRIMALITY

Florentin Smarandache
University of New Mexico
200 College Road
Gallup, NM 87301, USA
E-mail: smarand@unm.edu

**Abstract**: In this article we present four necessary and sufficient conditions for a natural number to be prime.

**Theorem 1.** Let $p$ be a natural number, $p \geq 3$: $p$ is prime if and only if $(p-3)! \equiv \dfrac{p-1}{2} \pmod{p}$.

*Proof:*

*Necessity:* $p$ is prime $\Rightarrow (p-1)! \equiv -1 \pmod{p}$ conform to Wilson's theorem. It results that $(p-1)(p-2)(p-3)! \equiv -1 \pmod{p}$, or $2(p-3)! \equiv p-1 \pmod{p}$. But $p$ being a prime number $\geq 3$ it results that $(2, p) = 1$ and $\dfrac{p-1}{2} \in \mathbb{Z}$. It has sense the division of the congruence by 2, and therefore we obtain the conclusion.

*Sufficiency:* We multiply the congruence $(p-3)! \equiv \dfrac{p-1}{2} \pmod{p}$ with $(p-1)(p-2) \equiv 2 \pmod{p}$ (see [1], pp. 10-16) and it results that $(p-1)! \equiv -1 \pmod{p}$, from Wilson's theorem, which makes us conclude that $p$ is prime.

**Lemma 1.** Let $m$ be a natural number $> 4$. Then $m$ is a composite number if and only if $(m-1)! \equiv 0 \pmod{m}$.

*Proof:*

The sufficiency is evident conform to Wilson's theorem.

*Necessity:* $m$ can be written as $m = a_1^{\alpha_1} \ldots a_s^{\alpha_s}$, where $a_i$ are positive prime numbers, two by two distinct and $\alpha_i \in \mathbb{N}^*$, for any $i$, $1 \leq i \leq s$.

If $s \neq 1$ then $a_i^{\alpha_i} < m$, for any $i$, $1 \leq i \leq s$.

Therefore $a_1^{\alpha_1} \ldots a_s^{\alpha_s}$ are distinct factors in the product $(m-1)!$ thus $(m-1)! \equiv 0 \pmod{m}$.

If $s = 1$ then $m = a^\alpha$ with $\alpha \geq 2$ (because $m$ is non-prime). When $\alpha = 2$ we have $a < m$ and $2a < m$ because $m > 4$. It results that $a$ and $2a$ are different factors in $(m-1)!$ and therefore $(m-1)! \equiv 0 \pmod{m}$. When $\alpha > 2$, we have $a < m$ and $a^{\alpha-1} < m$, and $a$ and $a^{\alpha-1}$ are different factors in the product $(m-1)!$.

Therefore $(m-1)! \equiv 0 \pmod{m}$ and the lemma is proved for all cases.

**Theorem 2.** Let $p$ be a natural number greater than 4. Then $p$ is prime if and only if $(p-4)! \equiv (-1)^{\left[\frac{p}{3}\right]+1} \cdot \left[\frac{p+1}{6}\right] \pmod{p}$, where [x] is the integer part of x, i.e. the largest integer less than or equal to x.

*Proof:*

*Necessity:* $(p-4)!(p-3)(p-2)(p-1) \equiv -1 \pmod{p}$ from Wilson's theorem, or $6(p-4)! \equiv 1 \pmod{p}$; $p$ being prime and greater than 4, it results that $(6, p) = 1$.

It results that $p = 6k \pm 1$, $k \in \mathbb{N}^*$.

A) If $p = 6k - 1$, then $6 | (p+1)$ and $(6, p) = 1$, and dividing the congruence $6(p-4)! \equiv p+1 \pmod{p}$, which is equivalent with the initial one, by 6 we obtain:

$$(p-4)! \equiv \frac{p+1}{6} \equiv (-1)^{\left[\frac{p}{3}\right]+1} \cdot \left[\frac{p+1}{6}\right] \pmod{p}.$$

B) If $p = 6k + 1$, then $6 | (1-p)$ and $(6, p) = 1$, and dividing the congruence $6(p-4)! \equiv 1 - p \pmod{p}$, which is equivalent to the initial one, by 6 it results:

$$(p-4)! \equiv \frac{1-p}{6} \equiv -k \equiv (-1)^{\left[\frac{p}{3}\right]+1} \cdot \left[\frac{p+1}{6}\right] \pmod{p}.$$

*Sufficiency:* We must prove that $p$ is prime. First of all we'll show that $p \neq \mathcal{M}6$.

Let's suppose by absurd that $p = 6k$, $k \in \mathbb{N}^*$. By substituting in the congruence from hypothesis, it results $(6k-4)! \equiv -k \pmod{6k}$. From the inequality $6k - 5 \geq k$ for $k \in \mathbb{N}^*$, it results that $k | (6k-5)!$. From $2 | (6k-4)$, it results that $2k | (6k-5)!(6k-4)$. Therefore $2k | (6k-4)!$ and $2k | 6k$, it results (conform with the congruencies' property) (see [1], pp. 9-26) that $2k | (-k)$, which is not true; and therefore $p \neq \mathcal{M}6$.

From $(p-1)(p-2)(p-3) \equiv -6 \pmod{p}$ by multiplying it with the initial congruence it results that: $(p-1)! \equiv (-1)^{\left[\frac{p}{3}\right]} 6 \cdot \left[\frac{p+1}{6}\right] \pmod{p}$.

Let's consider lemma 1; for $p > 4$ we have:

$$(p-1)! \equiv \begin{cases} 0 \pmod{p}, & \text{if } p \text{ is not prime;} \\ -1 \pmod{p}, & \text{if } p \text{ is prime;} \end{cases}$$

a) If $p = 6k + 2 \Rightarrow (p-1)! \equiv 6k \not\equiv 0 \pmod{p}$.
b) If $p = 6k + 3 \Rightarrow (p-1)! \equiv -6k \not\equiv 0 \pmod{p}$.
c) If $p = 6k + 4 \Rightarrow (p-1)! \equiv -6k \not\equiv 0 \pmod{p}$.

Thus $p \neq \mathcal{M}6 + r$ with $r \in \{0, 2, 3, 4\}$.

It results that $p$ is of the form: $p = 6k \pm 1$, $k \in \mathbb{N}^*$ and then we have: $(p-1)! \equiv -1 \pmod{p}$, which means that $p$ is prime.

**Theorem 3.** If $p$ is a natural number $\geq 5$, then $p$ is prime if and only if $(p-5)! \equiv rh + \dfrac{r^2-1}{24} \pmod{p}$, where $h = \left[\dfrac{p}{24}\right]$ and $r = p - 24h$.

*Proof:*
*Necessity:* if $p$ is prime, it results that:
$(p-5)!(p-4)(p-3)(p-2)(p-1) \equiv -1 \pmod{p}$ or
$24(p-5)! \equiv -1 \pmod{p}$.
But $p$ could be written as $p = 24h + r$, with $r \in \{1, 5, 7, 11, 13, 17, 19, 23\}$, because it is prime. It can be easily verified that $\dfrac{r^2-1}{24} \in \mathbb{Z}$.
$$24(p-5)! \equiv -1 + r(24h+r) \equiv 24rh + r^2 - 1 \pmod{p}.$$
Because $(24, p) = 1$ and $24 \mid (r^2 - 1)$ we can divide the congruence by 24, obtaining: $(p-5)! \equiv rh + \dfrac{r^2-1}{24} \pmod{p}$.

*Sufficiency:* $p$ can be written $p = 24h + r$, $0 \leq r < 24$, $h \in \mathbb{N}$.
Multiplying the congruence $(p-4)(p-3)(p-2)(p-1) \equiv 24 \pmod{p}$ with the initial one, we obtain: $(p-1)! \equiv r(24h + r) - 1 \equiv -1 \pmod{p}$.

**Theorem 4.** Let's consider $p = (k-1)!h + 1$, $k > 2$ a natural number.
Then: $p$ is prime if and only if
$$(p-k)! \equiv (-1)^{h + \left[\frac{p}{h}\right]+1} \cdot h \pmod{p}.$$

*Proof:* $(p-1)! \equiv -1 \pmod{p} \Leftrightarrow (p-k)!(-1)^{k-1}(k-1)! \equiv -1 \pmod{p} \Leftrightarrow (p-k)!(k-1)! \equiv (-1)^k \pmod{p}$.
We have: $((k-1)!, p) = 1$ \qquad (1)
A) $p = (k-1)!h - 1$.
a) $k$ is an even number $\Rightarrow (p-k)!(k-1)! \equiv 1 + p \pmod{p}$, and because of the relation (1) and $(k-1)! \mid (1+p)$, by dividing with $(k-1)!$ we have: $(p-k)! \equiv h \pmod{p}$.
b) $k$ is an odd number $\Rightarrow (p-k)!(k-1)! \equiv -1 - p \pmod{p}$ and because of the relation (1) and $(k-1)! \mid (-1-p)$, by dividing with $(k-1)!$ we have: $(p-k)! \equiv -h \pmod{p}$
B) $p = (k-1)!h + 1$
a) $k$ is an even number $\Rightarrow (p-k)!(k-1)! \equiv 1 - p \pmod{p}$, and because $(k-1)! \mid (1-p)$ and of the relation (1), by dividing with $(k-1)!$ we have: $(p-k)! \equiv -h \pmod{p}$.
b) $k$ is an odd number $\Rightarrow (p-k)!(k-1)! \equiv -1 + p \pmod{p}$, and because $(k-1)! \mid (-1+p)$ and of the relation (1), by dividing with $(k-1)!$ we have $(p-k)! \equiv h \pmod{p}$.

Putting together all these cases, we obtain: if $p$ is prime, $p = (k-1)!h \pm 1$, with $k > 2$ and $h \in \mathbb{N}^*$, then

$$(p-k)! \equiv (-1)^{h + \left[\frac{p}{h}\right] + 1} \cdot h \pmod{p}.$$

*Sufficiency:* Multiplying the initial congruence by $(k-1)!$ it results that:

$$(p-k)!(k-1)! \equiv (k-1)!h \cdot (-1)^{\left[\frac{p}{h}\right]+1} \cdot (-1)^k \pmod{p}.$$

Analyzing separately each of these cases:
A) $p = (k-1)!h - 1$ and
B) $p = (k-1)!h + 1$, we obtain for both, the congruence:

$$(p-k)!(k-1)! \equiv (-1)^k \pmod{p}$$

which is equivalent (as we showed it at the beginning of this proof) with $(p-1)! \equiv -1 \pmod{p}$ and it results that $p$ is prime.

**REFERENCE:**

[1]   Constantin P. Popovici, "Teoria numerelor", Editura Didactică şi Pedagogică, Bucharest, 1973.

[Published in "Gazeta Matematică", Bucharest, Year XXXVI, No. 2, pp. 49-52, February 1981]